\documentclass{amsart}

\usepackage[all]{xy}
\usepackage{hyperref}
\hypersetup{
    colorlinks=true, 
    linkcolor=blue,    
    citecolor=blue,       
    filecolor=blue,     
    urlcolor=blue           
}
\usepackage[lite]{amsrefs}
\usepackage{amssymb}
\usepackage{stmaryrd}
\usepackage{mathrsfs}

\newcommand{\loc}{\operatorname{loc}}

\newcommand{\supp}{\operatorname{supp}}

\newcommand{\Spec}{\operatorname{Spec}}

\newcommand{\Hom}{\operatorname{Hom}}

\newcommand{\relint}{\operatorname{relint}}

\newcommand{\coker}{\operatorname{coker}}
\newcommand{\Pol}{\operatorname{Pol}}
\newcommand{\CaDiv}{\operatorname{CaDiv}}

\newcommand{\ord}{\operatorname{ord}}

\newcommand{\triv}{\operatorname{triv}}

\renewcommand{\div}{\operatorname{div}}

\newcommand{\pp}{\mathbb{P}}

\newcommand{\qq}{\mathbb{Q}}
\newcommand{\zz}{\mathbb{Z}}

\newcommand{\cc}{\mathbb{C}}

\newtheorem{introthm}{Theorem}
\newtheorem{introcor}[introthm]{Corollary}

\newtheorem{theorem}{Theorem}[section]
\newtheorem{lemma}[theorem]{Lemma}
\newtheorem{proposition}[theorem]{Proposition}

\theoremstyle{definition}

\newtheorem{definition}[theorem]{Definition}

\newtheorem{remark}[theorem]{Remark}

\theoremstyle{remark}

\numberwithin{equation}{section}

\begin{document}

\title{On the topology of rational
$\mathbb{T}$-varieties \protect\\ of complexity one}
\author[A.~Laface]{Antonio Laface}
\address{
Departamento de Matem\'atica,
Universidad de Concepci\'on,
Casilla 160-C,
Concepci\'on, Chile}
\email{alaface@udec.cl}

\author[A. ~Liendo]{Alvaro Liendo}
\address{
Instituto de Matem\'atica y F\'isica,
Universidad de Talca,
Casilla 721,
Talca, Chile}
\email{aliendo@inst-mat.utalca.cl}

\author[J.~Moraga]{Joaqu\'in Moraga}
\address{
Department of Mathematics, University of Utah, 155 S 1400 E, 
Salt Lake City, UT 84112}
\email{moraga@math.utah.edu}

\thanks{
  \mbox{\hspace{11pt}}{\it Key words}: T-varieties, Deligne-Hodge
  polynomials, Chow rings, torus actions, topology of T-varieties.\\
  \mbox{\hspace{11pt}} The first author was partially supported by
  Proyecto FONDECYT Regular N. 1150732 and project Anillo ACT 1415 PIA
  Conicyt. The second author was partially supported by Proyecto
  FONDECYT regular N. 1160864.}

\maketitle

\medskip 

\begin{abstract}
We generalize classical results about 
the topology of toric varieties to the case of 
projective $\qq$-factorial $\mathbb{T}$-varieties 
of complexity one using the language of divisorial fans.  
We describe the Hodge-Deligne polynomial in
the smooth case, the cohomology ring and the 
Chow ring in the contraction-free case.
\end{abstract}

\setcounter{tocdepth}{1}
\tableofcontents

\section*{Introduction}
We study the topology
of {\em complex $\mathbb T$-varieties of 
complexity one}: normal complex
varieties equipped with an effective action 
of a codimension one torus $\mathbb T
= (\cc^*)^{n-1}$, where $n$ is the dimension
of the variety.
Such varieties are a natural generalization 
of $\mathbb T$-varieties of complexity zero,
i.e. toric varieties, and, as the latter, admit
a combinatorial description by means of the
language of divisorial fans which naturally
generalizes the language of fan of cones
used for toric varieties.
This language has been fully developed
in~\cite{AH05,AHH} and we will make use 
of it in this paper (see Section~\ref{background}).
Given a $\mathbb T$-variety of complexity
one $X$ there is a commutative diagram
of equivariant maps
\[
 \xymatrix{
  \widetilde{X}\ar[dr]_-{\pi}\ar[rr]^-r 
  && X\ar@{-->}[dl]\\
  & Y
 }
\]
where the dashed arrow is a rational 
map which realizes the Chow quotient for 
the torus action, the variety $Y$
is a smooth curve and $r$ is an 
equivariant resolution of indeterminacy.
It is known that $Y$ is rational if and only if 
$X$ is rational, since $X$ contains an open
subset isomorphic to $U\times (\cc^*)^k$, 
where $U$ is an open subset of $Y$.
The general fiber of $\pi$ is a toric variety
and $\pi$ is locally toric, 
this means that each fiber admits 
a formal toric neighborhood~\cite{LS}*{Proposition 2.6},
and we will exploit this property to
prove some of our results.
Moreover, there is a maximal open
subset $U\subseteq Y$ over which
$\pi$ is a trivial bundle.
The variety $X$ is 
{\em contraction-free} 
if $\widetilde X = X$. 

The topology of toric varieties has been 
studied by several authors, see
~\cite{CLS,DJ,Franz,FS,Fu} for example.
In this paper, we prove three theorems
about complexity one $\mathbb T$-varieties
which generalize previously known results
on toric varieties.
To begin with the statement of the 
first theorem we briefly recall 
that the definition of the Hodge-Deligne 
polynomial. This is a map 
$E\colon \text{Var}_\cc\to\zz[u,v]$
which factors through the Grothendieck ring
of varieties and induces a homomorphism 
of rings $K_0(\text{Var}_\cc)\to\zz[u,v]$.
Moreover, when $X$ is smooth projective then 
\[
 E(X) = \sum_{p+q \geq 0}(-1)^{p+q}h^{p,q}(X)u^pv^q,
\]
where $h^{p,q}(X)$ are the Hodge numbers
of $X$. 

Going back to $\mathbb T$-varieties we 
say that an orbit $\mathbb O\subseteq
\widetilde X$ for the torus action is 
{\em not contracted} if $r^{-1}(r(\mathbb O))
= \mathbb O$ and it is {\em contracted}
otherwise.
Observe that contracted orbits are 
exactly those ones which are mapped by
$r$ to the indeterminacy locus of 
$X\dashrightarrow Y$.
We denote by $n_k$ the number 
of non-contracted $k$-dimensional 
orbits in a general fiber of $\pi$, 
by $c_k$ the number of contracted ones in
the same fiber, and by $n_{U,k}$ the number 
non-contracted $k$-dimensional orbits in 
$\pi^{-1}(Y\setminus U)$, where $U$ 
is the open subset of $Y$ over which $\pi$
is a trivial bundle. These numbers 
can be computed directly in terms of
a defining divisorial fan for $X$
as explained in Remark~\ref{ident}.

\begin{introthm}\label{hodgetate}
Let $X$ be a smooth projective
$\mathbb{T}$-variety of complexity 
one with Chow quotient $Y$ of genus
$g$. Then the Deligne-Hodge polynomial
of $X$ is the following:
\[
 E(X)
 = 
 \sum_{k=0}^{\dim X}
 \Big(n_{k}(uv+g(u+v)+ |Y\setminus U|) + c_k + n_{U,k}\Big)
 (uv-1)^{k}.
\]
In particular, if $g=0$ then the Hodge structure on
each cohomology group of $X$ is of Hodge-Tate type,
i.e. $h^{p,q}(X)=0$ for $p\neq q$.
\end{introthm}

If $X$ is a $\mathbb Q$-factorial
projective rational $\mathbb{T}$-variety of 
complexity one and $X'\to X$ is a 
$\mathbb T$-equivariant resolution of 
singularities (see ~\cite[Theorem 0.1]{AW97}), 
then, as a consequence 
of Theorem~\ref{hodgetate},
the pure Hodge structure 
on each cohomology group of $X'$ is of 
Hodge-Tate type. 
By Poincar\'e duality and 
projection formula, which holds by the
$\mathbb Q$-factoriality assumption,
the cohomology of $X$ injects into the 
cohomology of $X'$, so that the former
is of Hodge-Tate type as well. Thus 
we have the following.

\begin{introcor}
If $X$ is a $\qq$-factorial projective rational 
$\mathbb{T}$-variety of complexity one
then the pure Hodge structure on each 
cohomology group of $X$ is of Hodge-Tate type.
\end{introcor}

The second result is about the canonical 
map from the rational Chow groups 
to rational Borel-Moore homology 
groups is an isomorphism for these 
varieties.
This property was recently proved for 
$\qq$-factorial toric varieties
admitting a proper toric morphism to 
an affine variety with a torus invariant point
(see ~\cite[Theorem A]{CMM}).

\begin{introthm}\label{chow}
Let $X$ be a $\qq$-factorial projective
contraction-free rational $\mathbb{T}$-variety
of complexity one.
Then the canonical map from the Chow 
groups to Borel-Moore homology
\[
 \Phi_k\colon A_k(X)_\qq \rightarrow H_{2k}(X,\qq)
\]
is an isomorphism, and the Borel-Moore
homology groups $H_k(X,\mathbb{Q})$
vanish for $k$ odd.
\end{introthm}

Finally, we determine a presentation for 
the rational Chow ring of a contraction-free 
rational $\mathbb T$-variety $X$ of complexity one.
Denote by $A^*(X)$ the Chow ring of $X$
and by $A^*(X)_\qq $ the $\qq$-algebra 
$A^*(X)\otimes_\zz\qq$.
The general fiber $X_0$ of 
$\pi\colon X\to\pp^1$ is a toric 
variety. Each prime torus-invariant
divisor of $X_0$ defines an 
invariant divisor of $X_0\times U
\simeq \pi^{-1}(U)$ and thus an 
invariant divisor of $X$ by taking 
closure. We call such invariant divisors 
{\em horizontal}. Each prime component 
of a fiber of $\pi$ over $\pp^1\setminus U$
is a {\em vertical} divisor.
Both sets, of horizontal and vertical divisors,
are finite by construction. We denote their union
by $\mathfrak D(X)$.
Define the homomorphism
of $\qq$-algebras
\[
 \qq[x_D\, :\, D\in\mathfrak D(X)] \to A^*(X)_\qq
 \qquad
 x_D\mapsto [D].
\]
We prove that this
homomorphism is surjective and determine its
kernel when the divisorial fan which
describes $X$ is {\em shellable}, a 
combinatorial technical
holds if for example $X$ is projective
and $\qq$-factorial.
In this case, we show that 
the kernel is the ideal $\mathscr I(X)$ 
generated by the linear forms
with integer coefficients
$\sum n_Dx_D$ such that $\sum n_DD$ 
is principal together with the square-free monomials 
$\prod_{D\in I}x_D$ such that the 
intersection $\bigcap_{D\in I}D$
is empty, where $I$ is a subset
of $\mathfrak D(X)$.

\begin{introthm}\label{chowgroup}
Let $X$ be a $\qq$-factorial projective
contraction-free rational $\mathbb{T}$-variety
of complexity one with $\mathfrak D(X)$ 
defined as above.
Then the assignment $x_D\mapsto [D]$
induces an isomorphism of 
finite dimensional $\mathbb{Z}_{\geq 0}$-graded
$\mathbb{Q}$-algebras
\[
  \qq[x_D\, :\, D\in\mathfrak D(X)]/\mathscr I(X) 
 \simeq
 A^*(X)_\qq.
\]
\end{introthm}

The above description of the Chow group can be recovered in the smooth
case using the techniques introduced by Brion
in~\cite[Section~3]{Brion}. In the smooth case the main technique is
the Bialynicki-Birula decomposition of smooth projective
T-varieties. Our approach is different, we generalize the shellability
condition introduced by Fulton to study the cohomology ring of
projective $\mathbb{Q}$-factorial toric varieties~\cite{Fu98}, see
also~\cite[Theorem 12.4.4]{CLS}.  The obstruction to generalize
Theorem~\ref{chowgroup} to any rational $\mathbb T$-variety of
complexity one, possibly non contraction-free, is that the rational
Chow ring is no longer generated in degree one, that is by divisors,
as shown in Remark~\ref{quadric} below.

The paper is organized as follows: 
In Section~\ref{background} we
introduce the combinatorial description 
of $\mathbb{T}$-varieties by means of
divisorial fans.
Section ~\ref{hodgedeligne}
deals with the Hodge-Deligne polynomials 
of rational complexity one $\mathbb{T}$-varieties, 
which allows us to compute the cohomology groups
in the smooth projective case.
In Section~\ref{cohomologyring} we prove 
Theorem~\ref{chow} while in Section~\ref{t3}
we prove Theorem~\ref{chowgroup}.

\section{Combinatorial description of $\mathbb T$-varieties}
\label{background}

\subsection{$\mathbb{T}$-varieties via divisorial fans}

In this section we recall the combinatorial description of
$\mathbb{T}$-varieties due to Altmann, Hausen and
S\"u\ss~\cite{AH05,AHH}, see also \cite{AIPSV} for a survey on known
results about $\mathbb{T}$-varieties. 

Let $N$ be a finitely generated
free abelian group of rank $n$ and let $M:=\Hom(N,\zz)$ be the dual of $N$.  We
denote by $N_\qq:=N\otimes_\zz\qq$ and $M_\qq:=M\otimes_\zz\qq$ the
associated rational vector spaces.
We denote by $\mathbb{T}_N$ 
the torus $\Spec \cc[M]$,
or simply by $\mathbb{T}$ 
when $N$ is clear from the context.
Let $\sigma$ be a pointed
polyhedral cone in $N_\qq$,
i.e. the only vector subspace 
contained in $\sigma$ is $0$.
For every convex polyhedron $\Delta
\subseteq N_\qq$ one defines its {\em recession cone} as
\[
\sigma(\Delta):=\{v\in N_\qq \, :\,  v+\Delta \subseteq \Delta\}
\] 
and we say that $\Delta$ is a 
{\em $\sigma$-polyhedron}
whenever $\sigma(\Delta)=\sigma$. 
The set $\Pol(\sigma)$ of all 
$\sigma$-polyhedra in $N_\qq$
is a semigroup with respect to
the Minkowski sum $\Delta + \Delta'
= \{p+p'\, :\, p\in\Delta,\, p'\in\Delta'\}$.
The neutral element is the polyhedra
$\{0\}$ and we will also consider the
empty polyhedron $\emptyset$
in $\Pol(\sigma)$ with addition rule 
$\emptyset + \Delta := \emptyset$ 
for all $\Delta\in\Pol(\sigma)$. 

Given a normal variety $Y$ one denotes
as usual by $\CaDiv(Y)$ the group of Cartier 
$\qq$-divisors.
A polyhedral divisor on
$(Y,N)$ is a formal sum
\[
 \mathcal{D}
 :=
 \sum_D \Delta_D \cdot D,
\] 
where $D$ varies along all 
the prime Cartier divisors of $Y$ 
and $\Delta_D$ is a convex
$\sigma$-polyhedra in $N_\qq$
which for all but finitely many 
$D$ equals $\sigma$. 
We denote by $\sigma(\mathcal{D})
= \sigma$ the common recession cone
of the polyhedra $\Delta_D$ and 
we call it the {\em recession cone} of $\mathcal{D}$.
The 
{\em locus} of $\mathcal D$ is
\[
 \loc(\mathcal{D})
 :=
 Y \setminus \bigcup_{\Delta_D=\emptyset} D
\]
and we say that $\mathcal{D}$ has {\em complete locus}
if $\loc(\mathcal{D})=Y$ holds. 
The {\em support} of $\mathcal{D}$,
denoted by $\supp(\mathcal D)$, is 
the subset of $\loc(\mathcal D)$ consisting of the
union of the prime divisors $D$ whose
coefficient $\Delta_D$ is non-zero, that
it not equal to $\sigma$. The 
{\em trivial locus} of $\mathcal{D}$ is
the complement of $\supp(\mathcal D)$
in $\loc(D)$.
Given a polyhedral divisor $\mathcal{D}$ on
$(Y, N)$ with recession cone $\sigma$ we have 
an evaluation map 
\[
 \mathcal{D}\colon \sigma^\vee \rightarrow 
 \CaDiv(Y)_\qq, \quad u\mapsto
 \mathcal{D}(u)
 :=\sum \min_{v\in \Delta_D} \langle u,v \rangle D.
\]
Observe that for any $u\in \sigma^\vee$,
the divisor $\mathcal{D}(u)$ has support contained
in $\supp(\mathcal{D})$.

\begin{definition}\label{p-div}
A $\qq$-divisor is {\em semiample} if it admits
a base point free multiple and it is {\em big}
if some multiple admits a section with 
affine complement.
Let $\mathcal{D}$ be a polyhedral divisor on $(Y,N)$ with recession
cone $\sigma$. Then $\mathcal{D}$ is a {\em p-divisor} if
$\mathcal{D}(u)$ is a semiample divisor for every $u\in \sigma^\vee$
and $\mathcal{D}(u)$ is big for 
every $u\in \relint(\sigma^\vee)$.
\end{definition}

Affine $\mathbb T$-varieties are 
described in the following way.
Given a p-divisor $\mathcal{D}$ 
on $(Y, N)$ with recession cone $\sigma$, we
denote by $\widetilde X(\mathcal{D})$ 
the relative spectrum on $\loc(\mathcal{D})$
of the coherent sheaf of algebras
\[
 \mathcal{A}(\mathcal{D}) :=\bigoplus_{u\in
  \sigma^\vee \cap M}\mathcal{O}_Y(\mathcal{D}(u)).
\]
Then $\widetilde X(\mathcal{D})$ 
is a normal variety with a $\mathbb{T}$-action 
induced by the $M$-grading.
The good quotient by this torus 
action $\pi\colon \widetilde X(\mathcal{D}) \to
\loc(\mathcal{D})$ is induced by the inclusion of sheaves
$\mathcal{O}_Y\to\mathcal{A}(\mathcal{D})$.  
Passing to global sections one gets
a normal affine $\mathbb{T}$-variety
\[
 X(\mathcal{D})
 :=
 \Spec \Gamma(\loc (\mathcal{D}),\mathcal{A}(\mathcal{D}))
\]
together 
with a proper birational morphism
$r\colon \widetilde{X}(\mathcal{D}) \rightarrow X(\mathcal{D})$.
The main result in \cite{AH05} states 
that every normal affine
$\mathbb{T}$-variety arises in this way.
We say that $X(\mathcal D)$ is 
contraction free if $r$ is the identity
map.
Given two p-divisors $\mathcal{D}$ 
and $\mathcal{D}'$ on $(Y, N)$	we
write $\mathcal{D}'\subseteq \mathcal{D}$ 
if each coefficient of $\mathcal{D}'$
is contained in the corresponding
coefficient of $\mathcal{D}$.
In this case, the inclusion 
induces a map 
$X(\mathcal{D}') \rightarrow X(\mathcal{D})$
and we say that
$\mathcal{D}'$ is a {\em face} of $\mathcal{D}$ if such map is an open embedding.  The
{ \em intersection } $\mathcal{D}\cap \mathcal{D}'$ of $\mathcal{D}'$
and $\mathcal{D}$ is the polyhedral divisor 
$\sum_D (\Delta_D \cap\Delta'_D)\cdot D$.

\begin{definition}\label{df}
  A {\em divisorial fan} $\mathcal{S}$ is a finite set of p-divisors
  on $(Y,N)$ such that the intersection of any two p-divisors of
  $\mathcal{S}$ is a face of both and $\mathcal{S}$ is closed under
  taking intersections.
When $N$ is clear from the context we say that $\mathcal{S}$
is a {\em divisorial fan} on $Y$.
\end{definition}

We denote by $X (\mathcal{S})$ the scheme obtained by gluing the
affine $\mathbb{T}$-varieties $X (\mathcal{D})$ and $X(\mathcal{D}')$ along the open
subvarieties $X(\mathcal{D} \cap \mathcal{D}')$ for each
$\mathcal{D},\mathcal{D}' \in \mathcal{S}$. As in the case of toric
varieties, $X(\mathcal{S})$ turns out to be a variety and since the
gluing is $\mathbb{T}$-equivariant, $X (\mathcal{S})$ is a
$\mathbb{T}$-variety \cite[Section~4.4]{AIPSV} . 
Analogously, one defines the variety
$\widetilde{X}(\mathcal{S})$ by gluing
the $\mathbb{T}$-varieties $\widetilde{X}(\mathcal{D})$
and $\widetilde{X}(\mathcal{D}')$ 
along the open subvarieties $\widetilde{X}(\mathcal{D}\cap\mathcal{D}')$ 
for each $\mathcal{D},\mathcal{D}'\in \mathcal{S}$.
We have an induced proper birational morphism 
$r\colon \widetilde{X}(\mathcal{S})\rightarrow X(\mathcal{S})$.
The set
$\{ \sigma(\mathcal{D}) \, :\,  \mathcal{D} \in \mathcal{S} \}$ 
is the {\em recession fan} $\Sigma (\mathcal{S})$ of $\mathcal{S}$.  The
{\em locus} of $\mathcal{S}$ is the open set
\[
 \loc (\mathcal{S}):=
 \bigcup_{\mathcal{D}\in\mathcal{S}} \loc (\mathcal{D}),
\] 
the {\em support} of $\mathcal{S}$ is $\supp(\mathcal{S}) :=
\bigcup_{\mathcal{D}\in \mathcal{S}} \supp (\mathcal{D})$ and the {\em
  trivial locus} of $\mathcal{S}$ is the open set $\triv(\mathcal{S}):=
\bigcap_{\mathcal{D} \in \mathcal{S}} \triv(\mathcal{D})$. The main result
in \cite{AHH} states that all normal $\mathbb{T}$-varieties 
arise in this way.

\section{The Hodge-deligne polynomial}
\label{hodgedeligne}

In this section we give explicit formulas 
for the class of a complexity one $\mathbb T$-variety
in the Grothendieck ring $K_0({\rm Var}_\cc)$
of complex varieties.
This is the free algebra generated
by the classes of complex varieties 
modulo the relations 
\[
 [X] := [X\setminus Z] + [Z]
 \qquad
 [X\times_\cc Y] := [X]\cdot [Y],
\]
where $Z$ is a closed subvariety of $X$.
In order to do this we will introduce the concept of toric bouquet
and recall the face-orbit correspondence described in ~\cite[Theorem 10.1]{AH05}.
See ~\cite{NS} for a reference on the Grothendieck ring of varieties.

\proof[Proof of Theorem~\ref{hodgetate}]
We begin by providing a formula
for the class of $X$ in the Grothendieck
ring of varieties.
If $\widetilde X_0$ is a fiber of the trivial fibration 
$V\to U$ then
\[
 [\widetilde X] 
 = 
 [U]\cdot [\widetilde X_0] + [\pi^{-1}(Y\setminus U)]
\]
holds.
Since $\widetilde X_0$ is a toric variety its class
is the sum of the classes of all of its 
disjoint orbits and the same holds 
for the class of $\pi^{-1}(Y\setminus U)$,
being the latter a union of toric bouquet.
The class of a $k$-dimensional orbit is
$(\mathbb L-1)^{k}$.
This gives the following formula
\[
[\widetilde{X}]= 
 \sum_{k=0}^n
 \Big((c_k+n_k)[U] + c_{U,k} + n_{U,k})\Big)
 (\mathbb L-1)^{k}.
\]
The effect of applying the equivariant 
birational morphism $r\colon \widetilde X\to X$
is that of collapsing all the contracted orbits
which now we want to count just once. 
More precisely, if $E\subseteq\widetilde X$ 
is the exceptional locus, not necessarily
of codimension one, then $r$ induces
an isomorphism $\widetilde X\setminus E
\to X\setminus r(E)$, so that 
\[
 [X] = [\widetilde X] - [E] + [r(E)].
\]

We claim that 
\[
[r^{-1}(\mathbb{O})] = [Y]\cdot [\mathbb{O}]
\]
holds in the Grothendieck ring.
Indeed, by~\cite[Corollarire 5.1]{Ver76}
we can stratify $Y$ into finitely many
locally closed subsets $Y=\cup_{i=0}^r Y_i$
such that the corresponding morphisms
$r^{-1}(\mathbb{O})_i \rightarrow Y_i$
are topologically trivial with fiber $\mathbb{O}$.
Since $Y$ is a curve, this locally closed stratification
consists of an open set $Y_0$ and finitely many points $Y_i$ with $i\in \{1,\dots, r\}$.
Clearly the fiber over the open set is isomorphic to $\mathbb{O}$
since the general fiber of $r^{-1}(\mathbb{O}) \rightarrow Y$
is isomorphic to $\mathbb{O}$.
On the other hand the fiber over $p_i$ for $i\in \{1,\dots, r\}$
is isomorphic to $\mathbb{O}$ as well, 
otherwise $r^{-1}(\mathbb{O})$ is not separated.
Thus, we have that following inequalities
$ [r^{-1}(\mathbb{O})] 
= [Y]\cdot [\mathbb{O}]$.
If we let $k = \dim \mathbb O$, it follows that 
\[
 [r^{-1}(\mathbb{O})] 
 = [Y]\cdot [\mathbb{O}]
 = ([U] + |Y\setminus U|) \cdot (\mathbb L-1)^k
\]
in the Grothendieck ring.
Thus we get the following formula for
the class of $X$:
\[
 [X] 
 = 
 \sum_{k=0}^n
 \Big(n_{k}[U] + c_k + n_{U,k}\Big)
 (\mathbb L-1)^{k}.
\]
The statement now follows by applying the
Hodge-Deligne polynomial to both sides
of the above equation and recalling that
it is a homomorphism of rings 
$E\colon K_0(\text{Var}_\cc) \to \zz[u,v]$.
\qed

\subsection{Combinatorial description
of the orbits}
\label{combinatorial}
In this subsection we recall 
the combinatorial description for
contracted (non-contracted)
orbits of a $\mathbb T$-variety.
Given a $\sigma$-polyhedron $\Delta \subseteq N_\qq$,
we denote by $\mathcal{V}(\Delta)$ its set of vertices
and by $\mathcal N := \mathcal{N}(\Delta)$ its {\em normal fan} consisting
of the regions where the function $u\mapsto \min_{v\in \Delta}\langle u,v\rangle$ is linear.
The cones of $\mathcal{N}(\Delta)$ are in one-to-one dimension-reversing
correspondence with the faces $F\preceq \Delta$ via the bijection
\[
 F \mapsto \lambda(F) 
 :=
 \{ u \in M_\qq \, :\,  \langle F, u\rangle=\min \langle \Delta, u\rangle \}.
\]
Given a $\sigma$-polyhedron $\Delta$ one defines 
the {\em affine toric bouquet}
$X(\Delta):=\Spec(\mathbb{K}[\mathcal{N}])$,
where $\mathbb{K}[\mathcal{N}]
:= \bigoplus_{u\in \sigma^\vee \cap M} \mathbb{K}\chi^u$
as a $\mathbb{K}$-vector space and the multiplication is given by
\[
\chi^u \cdot \chi^{u'}:=
\left\{
	\begin{array}{cl}
		\chi^{u+u'}  & \text{if $u$, $u'$ belong to a common cone of $\mathcal{N}$}, \\
		0 &  \text{otherwise}.
	\end{array}
\right. 
\]

The ring $\mathbb K[\mathcal{N}]$
is not an integral domain, hence $X(\Delta)$
is not a variety, but it is a scheme, 
nevertheless since $\mathbb K[\mathcal{N}]$
is reduced, 
$X(\Delta)$ has a decomposition
into a finite union of irreducible invariant 
toric varieties
\[
X(\Delta):=\bigcup_{v\in \mathcal{V}(\Delta)}  X(\qq_{\geq 0} \cdot ( \Delta - v)).
\]
Observe that the orbits of $X(\Delta)$
are in one-to-one dimension-reversing
correspondence with the faces $F\preceq \Delta$. 
We denote by $\Delta(k)$ the set of faces
of codimension $k$ and, given a face $F\in \Delta(k)$,
we denote by $\mathbb{O}_F$ its corresponding
torus orbit of dimension $k$.

A {\em polyhedral complex} on $N_\qq$
is a finite set $\triangle$ of polyhedra of 
$N_\qq$ such that, if $\Delta\in \triangle$
all the faces of $\Delta$ are in $\triangle$,
and the intersection of any two elements
of $\triangle$ is a face of both.
Given a polyhedral complex $\triangle$
on $N_\qq$ with maximal polyhedra $\Delta_1, \dots, \Delta_r$,
with respect to the inclusion,
we denote by $X(\triangle)$ the scheme
obtained by gluing $X(\Delta_i)$ and $X(\Delta_j)$
along $X(\Delta_i \cap \Delta_j)$ for each $i$ and $j$.
Observe that $X(\triangle)$ is a scheme with an action of the big torus $\mathbb{T}$.

\begin{remark}
For any maximal polyhedron $\Delta$ of $\triangle$
we have a well-defined preorder $\preceq$ of its faces.
The set of all these preorders induce 
a preorder on $\triangle$ denoted 
by $\preceq$ as well. 
The elements of $\triangle$ will be called {\em faces} 
since any such element is a face of a maximal polyhedra in $\triangle$.
\end{remark}

As before, there is a one-to-one 
dimension-reversing correspondence
between the faces of a polyhedral complex
$\triangle$ and the orbits of $X(\triangle)$.
Let $\mathcal{S}$ be a divisorial fan 
on a curve $Y$ which defines the 
$\mathbb T$-varieties $\widetilde X$ and 
$X$.
Given a point $p\in Y$ we define the 
{\em slice} $\mathcal{S}_p$ to be the 
polyhedral complex in $N_\qq$
which defines the fiber $\pi^{-1}(p)$.
Observe that if $p\in Y$ is a general
point the fiber is a toric variety and 
$\mathcal{S}_p$ is its fan (recession fan).
Any orbit of $\widetilde{X}$ corresponds
to a face $F$ of some slice $\mathcal{S}_p$.
Given a p-divisor 
$\mathcal{D}=\sum_{p\in Y} \Delta_p \cdot p$ 
on a curve $Y$, its {\em degree} is the polyhedron
\[
 \deg(\mathcal{D})=\sum_{p\in Y} \Delta_p \subseteq N_\qq.
\]

\begin{definition}
A face $F$ of $\mathcal{S}_p$
is {\em contracted} if there exists 
a p-divisor $\mathcal{D}\in \mathcal{S}$
of complete locus
such that $\sigma(F)\subseteq \sigma(\mathcal{D})$ and
$\deg(\mathcal{D})\cap \sigma(F)\neq \emptyset$.
In this case, we say that $\sigma(F)$
is a {\em contracted cone}.
\end{definition}

\begin{remark}\label{ident}
Observe that the definition of contracted
face only depends on the recession cone of that
face. In other words any two faces having the 
same recession cone, and possibly living
in different slices, are either both contracted 
or both not contracted. By~\cite[Theorem 10.1]{AH05}
an orbit $\mathbb O_F$ is contracted, i.e. 
$\mathbb O_F\subsetneq r^{-1}(r(\mathbb O_F))$,
if and only if the face $F$ is contracted.
In particular the three numbers $n_k$, $c_k$ 
and $n_{U,k}$ are directly computed from
the divisorial fan $\mathcal S$ which defines
the $\mathbb T$-variety $X$.
\end{remark}

\section{The Cohomology ring}
\label{cohomologyring}
From now on we will assume $X$ to be a
complete, rational, $\qq$-factorial, contraction-free
$\mathbb T$-variety of complexity one.
All these conditions can be checked 
directly on a divisorial fan $\mathcal S$ 
which defines $X$, see for example 
~\cite[Example 2.5]{LS}.
By the assumptions the Chow quotient is 
$Y = \pp^1$ and the quotient by the 
torus action is a morphism
\[
 \pi\colon X\to \pp^1.
\]
We will denote by $U\subseteq \pp^1$ the
maximal open subset where $\pi$ is 
a trivial bundle.
We will also denote by $\mathcal S$ the divisorial
fan on $\pp^1$ which defines $X$.
Our aim here is to describe the rational
cohomology ring of $X$.
The results of this sections generalizes
those already known in the toric
framework~\cite[Chapter 5.]{Fu}.
We begin by studying the chow groups 
of a fiber of $\pi$. Recall that
any such fiber is a {\em toric bouquet}, that
is a scheme of the form $X(\triangle)$,
for some polyhedral complex $\triangle$.
If the fiber is general then the polyhedral
complex is a fan, the recession fan, and the
fiber is the corresponding toric variety.

\begin{proposition}\label{chowbouquet}
Let $X(\triangle)$ be a toric bouquet.
Then the Chow group $A_k(X(\triangle))$ 
is generated by the orbit closures
of dimension $k$.
\end{proposition}

\begin{proof}
Let $X(\triangle)_i$ be the union of all the orbits of 
$X(\triangle)$ of dimension at most $i$. Then $X(\triangle)_i$ 
is a closed subset of $X(\triangle)$ and 
$X(\triangle)$ admits the filtration
$X(\triangle)=X(\triangle)_n \supseteq X(\triangle)_{n-1} \supseteq 
\dots \supseteq X(\triangle)_{-1}=\emptyset$.
By~\cite[Proposition 1.8]{Fu98}
for each $i$ there is an exact sequence
\[
 \xymatrix@1{
  A_k(X(\triangle)_{i-1})\ar[r] 
  & A_k(X(\triangle)_i) \ar[r]
  & A_k( X(\triangle)_i\setminus X(\triangle)_{i-1})\ar[r]
  & 0.
  }
\]
Observe that the second
map admits a section as a map
of sets: the class of an $i$-dimensional 
orbit of $X(\triangle)$ is sent to the class of its 
closure. The statement follows by 
induction on $i$.
\end{proof}

Recall that each orbit $\mathbb O
\subseteq X$ corresponds to a unique
face $F$ of a slice $\mathcal S_p$ of 
the divisorial fan $\mathcal S$.
In case $p$ is a general point the face
$F$ is a cone of the recession fan $\Sigma$,
which describes the toric fibers of
$\pi\colon X\to \pp^1$. To stress
this correspondence we will write
$\mathbb O_F$ for the orbit and 
$V(F)$ for its closure.
An orbit closure $V(F)$ is a 
{\em vertical subvariety} if it 
is mapped to $Y\setminus U$ by $r$.
Its class is a {\em vertical cycle}.
If $F = v$ is a vertex we denote by 
$D_{(p,v)}$ the vertical divisor $V(v)$.
A {\em horizontal subvariety} is the 
closure of $U\times\mathbb O_F$
for an orbit $\mathbb O_F$ which is 
mapped to $U$ by $r$.
We denote it by $V_Y(F)$ and call
its class a {\em horizontal cycle}.
In this case, $F$ is a cone of the recession fan.
If $F$ is a one-dimensional cone,
then $V_Y(F)$ is an invariant prime
divisor and we will denote it by $D_F$.

\begin{proposition}
\label{celldec}
The Chow group $A_k(X)$ is generated 
by the horizontal and vertical cycles of 
dimension $k$.
\end{proposition}

\begin{proof}
Let $X_0$ be a general fiber of 
the quotient morphism $\pi$. Since
$\pi$ is trivial over the open subset
$U\subseteq\pp^1$ we have an equivariant
isomorphism $V\simeq U\times X_0$.
Thus the $k$-th Chow group of 
$V$ is generated by horizontal
cycles.
The complement $X\setminus V$
is a union of toric bouquet and so,
by Proposition~\ref{chowbouquet} 
its $k$-th Chow group is generated 
by vertical cycles of dimension $k$.
The statement then follows by the 
exact sequence
\[
 \xymatrix@1{
  A_k(X\setminus V)\ar[r] 
  & A_k(X) \ar[r]
  & A_k(V)\ar[r]
  & 0.
  }
\]
\end{proof}

\begin{definition}
\label{celldec}
A toric bouquet $X(\triangle)$ admits a {\em cellular
decomposition} if it admits a filtration by 
invariant closed subsets
\[
 X(\triangle) = Z_k\supseteq Z_{k-1}\supseteq
 \cdots\supseteq Z_{0}\supseteq Z_{-1} 
 = \emptyset
\]
such that each $Y_i := Z_i\setminus Z_{i-1}$,
is a quotient of an affine space by 
a finite abelian group.
\end{definition}

\begin{remark}
\label{t3}
Observe that in the definition of cellular 
decomposition the closed subsets $Z_i$ 
can be reducible so that the closure
$\overline Y_i$ can be strictly contained 
in $Z_i$. Our definition of cellular 
decomposition is slightly weaker than
the one given in~\cite{Fu98}*{Example 1.9.1}
where each $Y_i$ is isomorphic to an affine
space. The statement~\cite{Fu98}*{Example 19.1.11(b)}, 
which is the content of our Theorem~\ref{chow},
remains true also with our definition.
A proof of this fact follows  
the argument given in
~\cite{Fu}*{Pag. 103--104}.
\end{remark}

\begin{proposition}\label{isomshell}
Let $X(\triangle)$ be a toric bouquet which admits a
cellular decomposition.
Then the classes of the closure of the 
$k$-dimensional varieties $Y_i$ form a vector
space basis of the rational Chow group $A_k(X(\triangle))_\qq$.
\end{proposition}

\begin{proof}
The statement is proved with
a similar argument to the one given
in the proof of Proposition~\ref{chowbouquet}.
In particular we observe that 
each $k$-dimensional $Y_i$ is
isomorphic to a quotient $\mathbb A^{k}/G$,
for some finite abelian group $G$ so that
$A_k(Y_i)_\qq\simeq\qq$ by
~\cite{Fu98}*{Example 1.7.6}.
\end{proof}

We recall a formula from
~\cite{LS}*{Example 2.5}.
Let $V_Y(\sigma)$
be a horizontal invariant subvariety defined
by the cone $\sigma$ of the recession fan of $X$
and let $p\in\pp^1\setminus U$ be a point
in the complement of the trivial locus $U$ of 
the quotient map $\pi\colon X\to \pp^1$.
Then the following equality holds in the
group of cycles of $X$:
\begin{equation}
\label{speci}
 V_Y(\sigma) \cap \pi^{-1}(p)
 \, =\,  
 \sum_{\{F\in\mathcal{S}_p\, :\, \sigma(F) = \sigma\}} \mu(v(F)) V(F),
\end{equation}
where $\mu(v)$ is the smallest positive 
integer such that the vector $\mu(v)v$ has
integer entries. The formula makes sense
since each face $F$ has at most one 
vertex $v(F)$, being the variety $\qq$-factorial
and $\dim(F) = \dim(\sigma(F))$.
As a consequence of the above formula we 
define the following specialization maps
at the level of Chow groups:
\begin{equation}
\label{sp}
 \phi_p\colon
 A_k(X_0)\to A_k(X_p)
 \qquad
 [V(\sigma)]\mapsto
 \sum_{\{F\in\mathcal{S}_p\, :\, \sigma(F) = \sigma\}} \mu(v(F)) [V(F)].
\end{equation}

\begin{remark}
\label{smap}
The specialization map~\eqref{sp}
can be geometrically understood 
as follows. Consider the 
fiber square
\[
 \xymatrix{
  X_p\ar[r]\ar[d] & X\ar[d]^-\pi\\
  p\ar[r] & \pp^1
 }
\]
and recall that there is a surjection
of Chow groups $A_k(X)\to A_k(X\setminus X_p)$
whose kernel is the image of $A_k(X_p)$.
According to~\cite{Fu98}*{\S 20.3}
and~\cite{Fu98}*{Remark 6.2.1}
the pullback $A_k(X)\to A_k(X_p)$ 
descends to a homomorphism
(specialization map)
\[
 A_k(X\setminus X_p)\to A_k(X_p)
 \qquad
 [Z] \mapsto [\overline Z\cap X_p],
\]
where $\overline Z$ is the closure
of $Z$ in $X$. The same construction
applies if we substitute $X$ with the
preimage of any open subset of $\pp^1$
which contains the point $p$. In particular
if we take $U\cup\{p\}$, recall that
$\pi^{-1}(U) \simeq U\times X_0$
and that $A_k(X_0\times U)\simeq A_k(X_0)$,
the corresponding specialization map
\[
 A_k(X_0)\to A_k(X_p)
 \qquad
 [Z]\mapsto [\overline{Z\times U}\cap X_p]
\]
equals~\eqref{sp} by~\eqref{speci}.
\end{remark}

\begin{definition}\label{con}
We say that $X$ is {\em shellable} if the
map~\eqref{sp} is injective for any 
$p\in\pp^1\setminus U$ and each
fiber of $\pi\colon X\to\pp^1$ admits
a cellular decomposition.
\end{definition}

\begin{proof}[Proof of Theorem~\ref{chow}]
In the proof we omit the rational coefficients 
in order to abbreviate notation.
Given a point $p\in\pp^1$ denote by
$X_p$ the fiber of $\pi$ over $p$ by
$X_0$ the fiber over a general point of 
$U\subseteq\pp^1$ and by $V = \pi^{-1}(U)$
the preimage of $U$.
We have a commutative diagram of 
Chow groups and Borel-Moore homology
with rational coefficients
\[
 \xymatrix@1{
 & A_k(X\setminus V)
 \ar[r]^-{j_*}\ar[d]^-\simeq & 
 A_k(X)\ar[r]\ar[d]^-{\Phi_k} &  
 A_k(V)\ar[d]^-\simeq\ar[r] & 0 \\
& H_{2k}(X\setminus V)\ar[r]^-{j_*^{\rm hom}}
&  H_{2k}(X)\ar[r] 
& H_{2k}(V)\ar[r] & 0
 }
\]
where the zero in the bottom row 
are due to the fact that each $X_p$
has trivial homology in odd dimension
by Proposition~\ref{isomshell}.
The first vertical arrow is an isomorphism 
by Proposition~\ref{isomshell}.
Using that $A_k(U) \rightarrow H_{2k}(U)$
is an isomorphism for each $k$, 
K\"unneth formula for Borel-Moore 
homology~\cite[Section 2.6.19]{CG} 
and Proposition~\ref{isomshell} it
follows that the third vertical
arrow is an isomorphism.
By the four lemma $\Phi_k$ 
is surjective. 
The commutativity of the first square
shows that $\ker(j_*)$ is mapped to 
$\ker(j^{\rm hom}_*)$ by the vertical 
isomorphism. In particular 
\[
 \dim\ker(j_*) \leq \dim\ker(j_*^{\rm hom})
\]
with equality if and only if 
$A_k(X)$ and $H_{2k}(X)$ have the
same dimension, i.e. if $\Phi_k$ is an 
isomorphism. By the long exact sequence
of homology there is a surjective map 
$H_{2k+1}(V)\to\ker(j_*^{\rm hom})$. 
Moreover, since $V \simeq U\times X_0$,
we have isomorphisms
$H_{2k+1}(V) 
\simeq H_1(U)\otimes_\qq H_{2k}(X_0)
\simeq H_1(U)\otimes_\qq A_{k}(X_0),
$
where the first one is by K\"unneth
formula, while the second one is due
to the fact that $X_0$ is toric.
In particular $\dim(H_1(U)\otimes_\qq A_{k}(X_0))
\geq \dim\ker(j_*^{\rm hom})$.
To conclude we show that there is
an injective linear map
\[
 \alpha\colon H_1(U)\otimes_\qq A_k(X_0)
 \to A_k(X\setminus V)
\]
whose image is contained in $\ker(j_*)$,
so that $\dim\ker(j_*) \geq \dim\ker(j_*^{\rm hom})$.
Observe that $H_1(U)$ is isomorphic to
the hyperplane of $\qq^{|\pp^1\setminus U|}$ 
consisting of vectors whose sum of coordinates 
equals $0$.
Using this identification we define $\alpha$ as
\[
 \alpha(m\otimes [Z])
 =
 \sum_{p\in\pp^1\setminus U}m_p\phi_p([Z]), 
\]
where the maps $\phi_p$ are defined
in~\eqref{sp}.
Since $X$ is a $\qq$-factorial
contraction-free projective variety, by 
Proposition~\ref{shellable} it is shellable
so that each $\phi_p$ is injective and
thus $\alpha$ is injective as well.
We conclude by observing that
the image of $\alpha$ is contained 
in $\ker(j_*)$ since by Remark~\ref{smap}
we have
\[
 j_*(\alpha(m\otimes [Z])) = [\div(f)]\cap [\overline{Z\times U}],
\]
where $f$ is pullback of a rational function 
of $\pp^1$ whose divisor is $\sum_{p\in \pp^1\setminus U}m_pp$.

\end{proof}

\subsection{Shellable polyhedral complexes}
We conclude the section with a 
combinatorial criterion for a toric
bouquet to admit a cellular 
decomposition.
Let $\triangle$ be a polyhedral complex
and let $X(\triangle)$ be the corresponding
toric bouquet.
\begin{definition}\label{shell}
A polyhedral complex $\triangle$ is {\em shellable}
if there exists an order on its maximal polyhedra
$\Delta_1,\dots,\Delta_k$ such that for each 
$i\in\{1,\dots,k\}$ the following set 
\[
 \{F\preceq\Delta_i\, :\, F \text{ is not contained in}
 \cup_{j<i}\Delta_i\}
\]
has a unique minimal element $G_i$.
\end{definition}

Observe that the definition of shellable
polyhedral complex depends on the 
given order of its maximal polyhedra.
It can happen that permuting the 
labelling the polyhedral complex is no
longer shellable. For an example of this
phenomenon see~\cite[Pag. 101]{Fu}.
The definition in Fulton book is slightly
different from our, but in the complete
case both definitions agree.
An example of shellable 
polyhedral complex is the fan 
of cones of a projective 
and simplicial toric variety
~\cite[Section 5.2, Lemma]{Fu}.
If we denote by
\begin{equation}
\label{defi}
 i(F) := {\rm min}\{j\in\zz_{>0}\, :\, F\preceq\Delta_j\},
\end{equation}
the definition of shellable polyhedral complex
is equivalent to say that the for each
face $F$ the set 
$\{F\in\triangle\, :\, i(F) = i\}$
has a unique minimal element $G_i$.

\begin{lemma}
\label{index}
Let $\triangle$ be a shellable polyhedral
complex with maximal polyhedra
$\Delta_1,\dots,\Delta_k$.
Then the following hold.
\begin{enumerate}
\item
$i(G_j) = j$.
\item
If $F\in \triangle$ then
$j = i(F)$ is the unique integer
such that $G_j\preceq F\preceq \Delta_j$.
\end{enumerate}
\end{lemma}
\begin{proof}
The first property is a direct consequence 
of the definition of $i$.
We now prove the second property.
Clearly $G_{i(F)}\preceq F\preceq \Delta_{i(F)}$ 
by the definition of $i$. This proves the 
existence part. To prove the unicity let
$j$ be such that $G_j\preceq F\preceq \Delta_j$,
then $j = i(G_j) \leq i(F) \leq j$.
\end{proof}

We say that a polyhedral complex $\triangle$
is {\em simplicial} if for each vertex $v\in \triangle$ 
and each maximal polyhedron $\Delta \in \triangle$ 
which contains $v$ the cone 
$\qq_{\geq 0}\cdot (v- \Delta)$ is simplicial.

\begin{proposition}
Let $\triangle$ be a shellable simplicial polyhedral
complex. Then the toric bouquet $X(\triangle)$
admits a cellular decomposition.
\end{proposition}
\begin{proof}
Let $\Delta_1, \dots , \Delta_k$ be the 
maximal polyhedra of $\triangle$.
For $1\leq i \leq k$ we define the 
invariant subvarieties of $X(\triangle)$ 
\[
 Y_i := \bigcup_{G_i \preceq F \preceq \Delta_i} \mathbb{O}_F,
 \qquad
 \qquad
 Z_i := Y_i \cup Y_{i+1} \cup \dots \cup Y_k.
\]
Observe that $Y_i = V(G_i) \cap U_{\Delta_i}$,
where $U_{\Delta_i}$ is the affine open
subset $\cup_{F \preceq \Delta_i} \mathbb{O}_F$.
The given definition coincides with the standard
definition of shellability of fans given in ~\cite{Fu}
for complete toric varieties.
Moreover, if a toric bouquet is shellable
then the order on its maximal polyhedra
induces an order on the cones
of the defining fan of each irreducible 
component of the bouquet.
Indeed the irreducible components of a 
toric bouquet are in bijection with the vertices
of the polyhedral complex. Putting the
origin of coordinates at a vertex and taking
the cone over the adjacent polyhedra
one gets the fan of cones which define
the corresponding irreducible component.
In particular an ordering of the maximal
polyhedra of the polyhedral complex
induces an ordering of the maximal cones
of each such component (see Subsection
~\ref{combinatorial}).
To prove that $Z_i$ is closed we will 
prove that it is equal to the finite union 
$\cup_{j>i}V(G_j)$
of closed subsets. One inclusion is
obvious since $Y_i\subseteq V(G_i)$ 
for each $i$. To prove the opposite 
inclusion let $p\in V(G_i)\setminus 
U_{\Delta_i}$, so that $p\notin Y_i$.
Then $p$ belongs to some orbit
$\mathbb O_F$ such that $G_i\preceq F$ 
and $F\not\preceq\Delta_i$.
By Lemma~\ref{index} we deduce that
$j := i(F)$ cannot be equal to $i$ and in 
particular $j > i$. 
Observe that 
\[
 G_j\preceq F\preceq \Delta_j
 \quad
 \text{if and only if}
 \quad
 \mathbb O_F\subseteq V(G_j) \cap U_{\Delta_j},
\]
so that the former inclusion holds.
Since $V(G_j) \cap U_{\Delta_j} = Y_j$ 
we conclude $p\in Y_j \subseteq Z_i$,
which proves the opposite inclusion.
It remains to show that if $\triangle$
is simplicial then $Y_i$ is the quotient
of an affine space by a finite group.
The open subset $U_{\Delta_i}$ of 
the toric bouquet $Z$ is an affine toric
variety, which is $\qq$-factorial due to
the simpliciality assumption.
Thus $Y_i$ is also toric, affine and 
$\qq$-factorial and the statement follows.
\end{proof}

\begin{proposition}
\label{shellable}
Let $X$ be a $\qq$-factorial contraction-free 
rational projective $\mathbb T$-variety of 
complexity one. Then $X$ is shellable.
\end{proposition}

\begin{proof}
Since $X$ is projective, the general fiber
$X_0$ of $\pi\colon X\to \pp^1$ is a projective
toric variety. Let $Q\subseteq M_\qq$ be its
defining compact polytope.
By~\cite[Section 14.4]{AIPSV}
every fibre $X_p$ of $\pi$ is defined
by a polyhedral subdivision of
$Q$, which we call $Q_p$. 
Let $Q'$ be a simplicial common refinement
of all the polyhedral subdivisions $Q_p$ for
$p\in\pp^1\setminus U$.
The scheme $X(Q')$ is a 
projective $\mathbb{Q}$-factorial toric bouquet,
meaning that every irreducible component is $\qq$-factorial.
Now, take an element $v\in N$ such that 
\[
 \langle v,u \rangle \neq \langle v,u' \rangle
 \qquad
 \text{for each pair of distinct vertices $u$ and $u'$ of $Q'$.}
\]
Let $\mathcal{S}_p$ be the polyhedral complex 
which defines the slice over $p$.
Observe that the set of maximal polyhedra
of $\mathcal{S}_p$ is in bijection with 
the set of vertices of $Q_p$.
Given a maximal polyhedron $\Delta$ of 
$\mathcal{S}_p$ we will denote by 
$u(\Delta)$ the corresponding vertex of
$Q_p$ via the above bijection.
Now, as $\Delta$ varies over the maximal 
polyhedra of $\mathcal{S}_p$, the values 
$\langle v,u(\Delta)\rangle$ are all distinct
and allow us to fix an order on the polyhedra.
By~\cite[Lemma, page 101]{Fu}
we induced a shellability on each 
polyhedral complex $\mathcal{S}_p$.
Recall that for each face $F$ of $\mathcal{S}_p$ 
we denote by $i(F)$ the minimum positive integer 
such that $F\preceq \Delta_{i(F)}$.
Thus, given two faces $G^0_n$ and $G^0_m$ 
of the fan $\mathcal S_0$ of the general 
fiber, with $n<m$, by construction 
\[
\min \{ i(F) \, :\,  F \in \mathcal{S}_p(G_n^0)\}
<
\min \{i(F) \, :\,  F \in \mathcal{S}_p(G_m^0)\}.
\]
This means that the representative matrix 
of $\phi_{p}$, with respect to the basis $B_0$
of $A_*(X_0)$, given in Proposition~\ref{isomshell},
and the basis $B_p$ of $A_*(X_p)$, contains 
a triangular submatrix of full rank,
so that $\phi_{p}$ is injective for each $p$.
\end{proof}

\section{The Chow ring}
\label{t3}

In this section we assume $X$ to be a 
complete $\qq$-factorial contraction-free 
rational $\mathbb T$-variety of complexity 
one. We will denote by $\pi\colon X\to Y$
the quotient morphism and by $\mathcal S$
the divisor fan of $Y$ which defines $X$.

\begin{lemma}
\label{intersection}
The rational Chow ring of $X$ is generated 
by classes of invariant divisors.
\end{lemma}

\begin{proof}
By~\cite{LS}*{Example 2.5} the 
variety $X$ is toroidal, that is
each point $p\in \pp^1$ admits an
analytic neighborhood $V\subseteq \pp^1$
such that $\pi^{-1}(V)$ is a toric variety.
Thus a vertical subvariety is intersection 
of invariant divisors since this is the case 
in a toric variety.
Consider now a horizontal subvariety
$V_Y(\sigma)$, which by definition is the 
closure of $U\times \mathbb O_\sigma$
for some cone $\sigma$ of the recession fan.
To conclude it suffices to prove the equality
\[
 V_Y(\sigma) = \bigcap_{\rho\in\sigma(1)} D_\rho,
\]
where $D_\rho$ is the horizontal divisor
defined by the one-dimensional subcone
$\rho$ of $\sigma$.
The inclusion $\subseteq$ is obvious. 
The opposite inclusion clearly holds in 
$V$, so it remains to prove it 
at the fibers over the points 
$p\in \pp^1\setminus U$.
This is an immediate consequence
of the fact that the fiber over
$p$ admits a $\qq$-factorial
toric neighborhood~\cite{LS}*{Proposition 2.6}
and thus, up to rational multiples, each cycle 
in it is a complete intersection.
\end{proof}

Before proceeding with the next lemmas
we briefly recall the formula for the divisor
of a torus-invariant rational function.
Such function is a product of thepullbacks
of a rational function $f\in\cc(\pp^1)$,
together with a character $\chi^u$ 
of the torus $T$, via the projection
$U\times T\to T$. For simplicity of
notation we will denote the function 
just by $f\chi^u$, omitting the pullbacks.
By~\cite{AIPSV}*{Theorem 26} we have
\begin{equation}
\label{principal}
 \div(f\chi^u)
 = 
 \sum_{\rho \in \mathcal R(\mathcal S)} \langle \rho, u \rangle D_\rho + 
 \sum_{v\in \mathcal V(\mathcal S)} 
 \mu (v)\left(\langle v, u \rangle + \ord_p(f)\right)  D_v.
\end{equation}
Observe that the first sum consists 
entirely of horizontal divisors, while 
the second sum contains the vertical
part of the principal divisor.

In what follows we will denote by 
$x_\rho$ the variable $x_{D_\rho}$ 
and by $x_v$ the variable $x_{D_v}$.
We also denote by $\mathcal R(\mathcal S)$
the set of primitive generators of the 
one-dimensional cones of the recession fan
of $\Sigma(S)$ and by $\mathcal V(\mathcal S)$
the set of vertices of $\mathcal S$.
With this notation the polynomial ring 
$\qq[x_D\, :\, D\in\mathfrak D(X)]$
becomes 
$\qq[x_w\, :\, w\in\mathcal R(\mathcal S)\cup
\mathcal V(\mathcal S)]$.
Given an irreducible invariant subvariety 
$Z$ of $X$ defined by a face $F$ of the
divisorial fan $\mathcal S$ we define 
the monomial
\[
 m_F
 = 
 \prod_{\rho\in\sigma(F)}x_\rho
 \cdot
 \prod_{v\in\mathcal V(F)}x_v.
\]
We denote by $\mathscr I(X)$ the ideal
of the polynomial ring $\qq[x_D\, :\, D\in\mathfrak D(X)]$
generated by all the above monomials
together with the linear polynomials
coming from the linear relations from 
the divisors of $\mathfrak D(X)$.

\begin{lemma}
\label{linear}
Let $w$ be either a vertex or a ray of a face
$F$ of the divisorial fan $\mathcal S$. 
Then there is a linear relation with 
rational coefficients
\[
 x_w
 \equiv 
 \sum_{v\in\mathcal V(\mathcal S)\setminus\mathcal V(F)} \alpha_v\, x_v
 +
 \sum_{\rho\in\mathcal R(\mathcal S)\setminus\mathcal R(F)}\alpha_\rho\, x_\rho
 \mod \mathscr I(X).
\]
\end{lemma}
\begin{proof}
First of all we recall the local description 
of the variety $\pi^{-1}(U_p)$, where $U_p$ 
is an euclidean neighbourhood of a point
$p\in\pp^1$. The variety is toric and 
its fan can be described as follows.
Let $\Sigma(\mathcal S)$
be the recession fan and $\mathcal S_p$ be the
slice over the point $p$. Form the 
rational vector space $N_\qq\oplus\qq$,
embed the recession fan in $N_\qq\oplus\{0\}$
and the slice in $N_\qq\oplus\{1\}$. 
The set of cones of the embedded recession fan
together with the cones on the faces of the 
embedded slice form the fan of $\pi^{-1}(U_p)$. 
The variety is $\qq$-factorial if and only if for
any face $F$ of $\mathcal S_p$ the following
cone
\[
 \tilde\sigma_F 
 =
 {\rm cone}(
 \sigma(F)\times\{0\}
 \cup
 \mathcal V(F)\times\{1\}
  )
 \subseteq N_\qq\oplus\qq
\]
is simplicial.
Since the variety $X$ is $\qq$-factorial
one can find
an element $u'\in M\oplus\zz$ which
vanishes on all the extremal rays
of $\tilde \sigma_F$ but $(w,1)$. 
Denote by $u$ the restriction of $u'$ to $M$.
Let $f$ be a rational function of $\pp^1$,
having zeroes and poles in $\pp^1\setminus U$,
and such that $\ord_{p}(f) = u'(0,1)$.
Observe that if $v$ is a vertex of the 
slice $\mathcal S_p$ containing $F$ 
then the expression
\[
 \langle v, u \rangle + \ord_p(f)
 =
 u(v) + \ord_p(f)
 =
 u'(v,0) + u'(0,1)
 =
 u'(v,1)
\]
vanishes exactly when $v\neq w$.
Then equation~\eqref{principal} and
our hypothesis on $u$ imply that
the linear polynomial
\[
 \alpha_w\, x_w +
 \sum_{v\in \mathcal V(\mathcal S)\setminus\mathcal V(F)} 
 \mu (v)\left(\langle v, u \rangle + \ord_p(f)\right) x_v +
 \sum_{\rho \in \mathcal R(\mathcal S)\setminus\mathcal R(F)}
 \langle \rho, u \rangle x_\rho
\]
is in $\mathscr I(X)$, where 
$\alpha_w = \mu(w)u'(w,1)\neq 0$. 
The statement follows.
\end{proof}

\begin{lemma}
\label{corresp}
Each monomial of $\qq[x_w\, :\, w\in
\mathcal R(\mathcal S)\cup
\mathcal V(\mathcal S)]$
is equivalent, modulo $\mathscr I(X)$,
to a linear combination of square-free 
monomials with rational coefficients.
Moreover a square-free monomial
not contained in $\mathscr I(X)$
is of the form $m_F$ for some face 
$F$ of $\mathcal S$.
\end{lemma}

\begin{proof}
Let $m = \prod_{w\in S}x_w^{a_w}$ be a 
monomial, where $S$ is a subset of 
$\mathcal V(\mathcal S)\cup \mathcal R(\mathcal S)$.
The proof of the first statement is by induction
on $\mu(m) = \sum_{w\in S}(a_w-1)$.
If $\mu(w) = 0$ then $m$ is square-free
and there is nothing to prove.
Otherwise either $m$ is in $\mathscr I(X)$
or $S$ is contained in $\mathcal V(F)\cup \mathcal R(F)$
for some face $F$ of the divisorial fan
$\mathcal S$. In the second case
let $w\in S$ be such that $a_w > 1$.
Applying Lemma~\ref{linear} to $x_w$ 
and multiplying both sides of the equation
by $m/x_w$ we get
\[
 m
 \equiv 
 \sum_{v\in\mathcal V(\mathcal S)\setminus\mathcal V(F)} \alpha_v\, x_v\frac{m}{x_w}
 +
 \sum_{\rho\in\mathcal R(\mathcal S)\setminus\mathcal R(F)}\alpha_\rho\, x_\rho\frac{m}{x_w}
 \mod \mathscr I(X).
\]
Since $S\subseteq \mathcal V(F)\cup \mathcal R(F)$
we deduce $\mu(m')<\mu(m)$ for each
monomial $m'$ on the right hand side of
the above expression. This proves
the first statement.

For the second statement observe that,
being $X$ a locally toric 
$\qq$-factorial variety, the invariant
subvariety represented by a square-free 
monomial not contained in $\mathscr I(X)$
corresponds to a cone $\tilde \sigma_F$
defined by a face $F$ as in the proof
of Lemma~\ref{linear}.
\end{proof}

\begin{lemma}\label{algmov}
Let $\emptyset\neq F_1\prec H\preceq F_2$
be three faces of the same slice $\mathcal S_p$.
Then 
there exist rational numbers $\alpha_F$ 
such that
\[
 m_{H} \equiv \sum_{F_1\preceq F\not\preceq F_2} \alpha_F\, m_F
 \mod \mathscr I(X),
\]
where the sum varies over all the faces
containing $F_1$ and not contained
in $F_2$.
\begin{proof}
By hypothesis either $H$ contains a vertex 
which is not in $F_1$ or $\sigma(H)$
contains a ray which is not a ray of 
$\sigma(F_1)$.
Assume we are in the 
first case and let $w\in\mathcal V(H)
\setminus\mathcal V(F_1)$.
Using the fact that $w\in \mathcal V(F_2)$
and multiplying both sides of the relation
in Lemma~\ref{linear} by $m_{H}/x_w$
we get the following expression

\[
 m_{H}
 \equiv 
 \sum_{v\in\mathcal V(\mathcal S)\setminus\mathcal V(F_2)} \alpha_v\, x_v\frac{m_{H}}{x_w}
 +
 \sum_{\rho\in\mathcal R(\mathcal S)\setminus\mathcal R(F_2)}\alpha_\rho\, x_\rho\frac{m_{H}}{x_w}
 \mod \mathscr I(X).
\]
By Lemma~\ref{corresp} each non-zero 
monomial is of the form $m_F$ for some 
face $F$ of $\mathcal S$.
Observe that $F_1\preceq F$
since $F$ contains all the vertices
and rays of $F_1$. Moreover 
$F$ is not a face of $F_2$
since it contains either a vertex or
a ray which is not in $F_2$.
A similar argument applies to the case
when $\sigma(H)$ contains a ray
which is not a ray of $\sigma(F_1)$.
The statement follows.
\end{proof}
\end{lemma}

To any $p\in\pp^1\setminus U$ we
associate an ideal 
$I_p  = \langle x_v +\mathscr I(X)\, :\, 
v\in\mathcal V(\mathcal S_p)\rangle$
of the quotient ring. Our next aim is to 
give a vector space basis of $I_p$.

\begin{lemma}\label{basislemma}
Assume that there is a cellular decomposition
of the toric bouquet $X_p = Z_k\supseteq Z_{k-1}\supseteq
 \cdots\supseteq Z_{0}\supseteq Z_{-1} 
 = \emptyset$, where the closure of 
 $Z_i\setminus Z_{i-1}$
corresponds to the face $G_i$ 
of the slice $\mathcal S_p$.
Then the ideal $I_p$ is generated
as a $\qq$-vector space by the monomials
$m_{G_1},\dots,m_{G_k}$.
\end{lemma}
\begin{proof}
By Lemma~\ref{linear} a vector space
basis of $I_p$ consists of square-free
monomials of the form $m_F$ for some
face $F$ of the divisorial fan $\mathcal S$.
Moreover, $F$ must be a face of the slice
$\mathcal S_p$ since the ideal generators
of $I_p$ are vertices and rays contained
in this slice.

Now let $F$ be a face of the slice $\mathcal S_p$.
By Lemma~\ref{index} there is a unique $i$
such that $G_i\preceq F\preceq \Delta_i$.
We prove, by descending induction on $i$,
that
\[
 m_F\equiv \sum_{j\geq i}\alpha_jm_{G_j^p}
 \mod \mathscr I(X), 
\]
with $\alpha_j\in\qq$ for any $j$. 
If $F=G_i^p$ we are done. Otherwise
the class of $m_{F}$ is equal to a class 
of the form $\sum_H \alpha_H\, m_H$
by Lemma~\ref{algmov} where each
face $H$ satisfies $G_j^p\preceq H\preceq \Delta_j$ 
for some $j>i$. We conclude by the induction
hypothesis.
\end{proof}

\begin{proof}[Proof of Theorem~\ref{chowgroup}]

Let $R$ be the quotient ring
$\qq[x_w\, :\, w\in\mathcal R(\mathcal S)\cup\mathcal V(\mathcal S)]
 /\mathscr I(X)$.
We have a well-defined surjective homomorphism
\[
 \phi \colon R \to A^*(X)_\qq
 \qquad
 x_w+\mathscr I(X) \mapsto [D_w],
\]
where $D_w$ is the invariant prime 
divisor defined by $w$.
We will produce a generating set of
the domain, as a vector space, whose 
image via $\phi$ is a vector space basis 
of the rational Chow ring.

Observe that the subring of $A^*(X)_\qq$
generated by $\{x_\rho \, :\,  \rho \in \mathcal R(\mathcal{S})\}$
is isomorphic to $A^*(X_0)_\qq$, 
and by Proposition~\cite[Chapter 5]{Fu} it is generated
as $\qq$-vector space by the classes of all 
the monomials $m_{G^0}$.
From now on we use the notation of 
the proof of Theorem~\ref{chow}.
The restriction of $\phi$ to the ideal 
$I_p\subseteq R$ induces a
surjective homomorphism
\[
 \alpha_p\colon I_p\rightarrow A^*(X_p)_\qq
\]
of $\qq$-vector spaces
which maps the set of classes of all 
monomials of the form $m_{G^p}$ to 
a basis of the codomain, by Proposition
~\ref{isomshell}.
Then $\alpha_p$
is an isomorphism by Lemma~\ref{basislemma}.
Consider a rational function $f\in \cc(\pp^1)$
of degree one with $\div(f) = p_0-p_i$.
Expanding the product $m_{F}\div(f)$ 
we deduce the following
\[
 m_{F}
 \equiv 
 \sum_{v\in \mathcal{V}(\mathcal{S}_p)} 
 m_{F}\mu(v)x_{v}
 \mod \mathscr I(X).
\]
For each $p$, we define the $\qq$-linear map
$\psi_p\colon I_0 \rightarrow I_p$ which sends the 
class of $m_F$ to the class of 
$\sum_{v\in \mathcal{V}(\mathcal{S}_p)}
 \mu(v) m_Fx_v$.
This makes the following diagram commutative
\[
 \xymatrix{
  I_0 \ar@{^{(}->}[r]^-{\psi_p}\ar[d]^-{\simeq} &
  I_p \ar@{^{(}->}[r]\ar[d]^-{\simeq} & 
  R \ar@{->>}[d]^-{\phi} \\ 
  A_*(X_0)_\qq \ar@{^{(}->}[r]^-{j_p}
  & A_*(X_p)_\qq \ar@{^{(}->}[r]
  & A_*(X)_\qq.
 }
\]
In particular for each $p\in\pp^1\setminus U$ 
we can choose a set of elements
$\mathcal{B}_p$ of $I_p$ which are linear combination
of the classes of the monomials $m_{G^p}$, such that
the image of $\mathcal{B}_p$ via $\phi$
is a basis of $\coker(j_p)$.
Observe that the union $\mathcal{B}_0\cup\mathcal{B}_p$,
where $\mathcal{B}_0$ consists of the
classes of the monomials $m_{G}$,
generates $I_p$ as a $\qq$-vector space for each $p$.
Being $X$ a $\qq$-factorial toroidal variety
by~\cite{LS}*{Example 2.5} we have that
$A^{n-k}(X)\simeq A_k(X)$
and using Theorem~\ref{chow} we have that
\begin{align*}
 A_k(X)_\qq
 & \simeq A_k(V)_\qq\oplus A_k(X\setminus V)_\qq/\ker(j_*)\\
 & \simeq A_{k-2}(X_0)\oplus A_k(X_0)
 \oplus\bigoplus_{p\in\pp^1\setminus U} \coker(j_p)_k,
\end{align*}
where we are making use of the isomorphism
$V\simeq U\times X_0$, of the fact that the 
Chow groups of $U$ are trivial, and 
where $\coker(j_p)_k$ denotes the $k$-graded
piece of $\coker(j_p)$.
We conclude that $\mathcal B_0
\cup\bigcup_{p\in\pp^1\setminus U} \mathcal{B}_p$
is a subset of $R$
which generates it as a $\qq$-vector space
and whose image is a $\qq$-basis of the Chow ring.
\end{proof}

\begin{remark}
\label{quadric}
Observe that the conclusion of Theorem
~\ref{chowgroup} is no longer true if we
substitute $\widetilde X(\mathcal{S})$ with
$X(\mathcal{S})$. As an example consider
the quadric $Q = V(x_1x_2+x_3x_4+x_5x_6)$
of $\pp^5$. It admits an effective action of 
$(\cc^*)^3$ and thus is a $\mathbb T$-variety of 
complexity one, so that $Q = X(\mathcal{S})$
for some divisorial fan $\mathcal S$ on 
$\pp^1$. On the other hand $Q$ is isomorphic
to the Pl\"ucker embedding of the Grassmannian
$G(2,4)$.
The dimension of the cohomology groups 
of $Q$ are well known (alternatively one can
easily compute them using Theorem~\ref{hodgetate}):
\[
h^0(Q)=1,\quad h^2(Q)=1,\quad h^4(Q)=2,\quad
h^6(Q)=1, \quad  h^8(Q)=1.
\]
Then, the Chow ring of $Q$ is generated in degree one and two.
\end{remark}

\end{document}